\documentclass[11pt]{amsart}
\usepackage{color}
\usepackage{amssymb,amsmath}
\usepackage{hyperref}
\usepackage{amsthm}
\usepackage{enumitem}
\usepackage{color}

\usepackage{wasysym}
\usepackage{graphicx}
\usepackage{arydshln}
 \usepackage{tikz}
 \usepackage{float}
 \usepackage{mathtools}

\newtheorem{theorem}{Theorem}[section]
\newtheorem{lemma}[theorem]{Lemma}
\newtheorem{example}[theorem]{Example}
\newtheorem{cor}[theorem]{Corollary}
\newtheorem{remark}[theorem]{Remark}

\newtheorem{claim}{Claim}

\newtheorem{definition}[theorem]{Definition}

\def\f{\noindent}
\def\h{\hfill  $\Box$\vspace{10pt}}
\def\p{\f {\bf Proof}\hskip10pt}

\def\tr{{\rm tr\;}}
\def\Ad{{\rm Ad}}
\def\Int{{\rm Int}}

\newcommand{\GL}{\mbox{\rm GL}}

\newcommand{\R}{\Bbb{R}}
\newcommand{\C}{\Bbb{C}}

\newcommand\norm[1]{\left\lVert#1\right\rVert}

\font\germ=eufm10

\def\g{{\mbox{\germ g}}}

\def\gl{{\mbox{\germ gl}}}

\def\h{{\mbox{\germ h}}}
\def\p{{\mbox{\germ p}}}
\def\k{{\mbox{\germ k}}}

\def\z{{\mbox{\germ z}}}

\begin{document}


\title{Curvature of matrix and reductive Lie groups}
\author{Luyining Gan}
\author{Ming Liao}
\author{Tin-Yau Tam}
\address{Department of Mathematics \& Statistics, Auburn University, Auburn, AL, 36849, USA. Email: \tt{lzg0027@auburn.edu}}
\address{Department of Mathematics \& Statistics, Auburn University, Auburn, AL, 36849, USA. Email: \tt{liaomin@auburn.edu}}
\address{Department of Mathematics \& Statistics, University of Nevada, Reno, NV 89557, USA. Email: \tt{ttam@unr.edu}}
\thanks{In Honor of Professor Jimmie Lawson: 50 Years at LSU}
\maketitle

\begin{abstract}
In this paper, we give a simple formula for sectional curvatures on the general linear group, which is also valid for many other matrix groups.
Similar formula is given for a reductive Lie group.
We also discuss the relation between commuting matrices and zero sectional curvature.

\bigskip
\f {\it Mathematics Subject Classification 2010:}  53B20, 14L35, 51N30

\f {\bf Keywords:} Curvature, general linear group, reductive Lie group, closed subgroup
\end{abstract}

The curvature provides important information about the geometric structure of a Riemannian manifold. 
For example, it is related to the rate at which two geodesics emitting from the same point move away from each other: the lower the curvature is, the faster they move apart (see Theorem~IX.5.1 in \cite[Chapter~IX.5]{Chavel1993}). Many important geometric and topological properties are implied by suitable curvature conditions.
However, the curvature is usually not easy to compute explicitly. In the case of a Lie group equipped with a left invariant Riemannian metric, Milnor \cite{Milnor} obtained an explicit formula for sectional curvatures,
but it is still quite complicated. To use it to compute a sectional curvature, one has to embed the spanning vectors of the section in an orthonormal frame and to compute the structure constants of the frame.
Although it simplifies in many special cases,
we have not seen a simple formula for the sectional curvature on the general linear group of matrices that is valid for all sections.

Our main result of this paper (Theorem~\ref{Th_general}) is a simple and direct formula for the sectional curvature on the general linear group equipped with the left invariant Riemannian metric induced by the Frobenius norm. 
This formula also holds on any matrix group that is invariant under transposition, such as the orthogonal group and the Lorentz group, as they are totally geodesic submanifolds of the general linear group.
Indeed, similar formula is valid for reductive Lie groups. More details will be given later.

Our paper is organized as follows. After the preliminary material is introduced in Section 1,  we establish our main result on the sectional curvature for the general linear group in Section 2.
In Section 3 we 
study the sectional curvature on $\GL(n,\R)$ when the two tangent vectors are commuting matrices in ${\gl}(n,\R)$.
In Section 4 an extension in the context of a  reductive Lie group $G$ is given. The results in Section 2 are then particular cases. The intent of the separation of Section 2 from Section 4 is to make the results more accessible to matrix and application oriented readers who may not be  familiar with reductive Lie group.

\section{Preliminary}

The set $M(n,\R)$ of $n\times n$ real matrices may be identified with the Euclidean space $\R^{n^2}$. The general linear group $G= \GL(n, \R)$ is the open subset containing only nonsingular matrices and is a Lie group under the matrix multiplication. Its tangent space $T_eG$ at the identity element $e$, which is just the identity matrix, is the Lie algebra of $G$ and is denoted by ${\g} = {\gl}(n, \R)$.
Note that ${\g}$ may be identified with $M(n,\R)$ and it is equipped with the Lie bracket $[u,v]=uv-vu$ satisfying the Jacobi identity
\begin{equation}\label{Eq_Jacobi}
[u, [v, w]]+[v, [w, u]]+[w, [u, v]] = 0.
\end{equation}
Let ${\g}$ be equipped with the Euclidean inner product, known as the Frobenius inner product, defined
by $\langle X, Y \rangle = \mathrm{tr} (X^\top Y)$, where $\mathrm{tr}(\cdot)$ stands for the trace and the superscript $\top$ denotes the transpose.

As a manifold, $G$ is covered by a single coordinate neighborhood via the exponential map $\exp: \g\to G$ and any $g\in G$ has coordinates $g_{ij}$, the matrix elements of $g$. For $h\in G$, let $L_h$: $G\to G$ denote the left translation $g\mapsto hg$.
A (smooth) vector field $X$ on $G$ is called left invariant if $DL_h(X(g))=X(hg)$ for any $g,h\in G$, where $DL_h$ is the differential map of $L_g$. A left invariant vector field $X$ on $G$ is completely determined by its value $X(e)$
at $e$.

For simplicity, given any $u\in{\g}$, we may use $u$ also for the unique left invariant vector field $X$ with $X(e)=u$. Note that in terms of local coordinates $g_{ij}$, the value of $u$ at $g$,
as a left invariant vector field, is $gu$. Moreover, the Lie bracket $[u,v]=uv-vu$ may be understood in terms of either matrix multiplication or vector field operation.

A Riemannian metric on $G$ is a smooth distribution of inner products $\langle\cdot,\cdot\rangle_g$ on tangent spaces $T_gG$ for $g\in G$. It is called left invariant if for any $g,h\in G$ and $x,y\in T_gG$,
$\langle DL_hx,DL_hy\rangle_{hg}=\langle x,y\rangle_g$. A left invariant Riemannian metric on $G$ is completely determined by the inner product $\langle\cdot,\cdot\rangle=\langle\cdot,\cdot\rangle_e$ at $e$.
By (5.3) in \cite{Milnor}, the covariant derivative $\nabla_uv$ under a left invariant metric may be determined by
\begin{equation} \label{covderiv}
\langle\nabla_uv,w\rangle = \frac{1}{2}(\langle [u,v],w\rangle - \langle [v,w],u\rangle - \langle [u,w],v\rangle).
\end{equation}

In the rest of the paper, we will assume that $G$ is equipped with the left invariant Riemannian metric determined by the Frobenius inner product $\langle\cdot,\cdot\rangle$ at ${\g}$.
Let $\norm{u}=\langle u,u\rangle^{1/2}$ be the associated norm.
It is easy to show $\langle [u,w],v\rangle=\langle w,[u^\top,v]\rangle$. Then by (\ref{covderiv}),
\[\langle\nabla_uv,w\rangle = \frac{1}{2}(\langle [u,v],w\rangle - \langle [v^\top,u],w\rangle - \langle [u^\top,v],w\rangle).\]
It follows that
\begin{equation}\label{Eq_covder}
	\nabla_uv = \frac{1}{2}([u,v] + [u,v^\top]+[v,u^\top]).
\end{equation}

Let $\mathcal{S}$ and $\mathcal{A}$ be the spaces of symmetric and skew-symmetric matrices in the Lie algebra ${\g}$ of $G$, respectively. 
Note that ${\g}=\mathcal{S}\oplus\mathcal{A}$
is an orthogonal direct sum. For any two subspaces ${\g}_1$ and ${\g}_2$ of ${\g}$, let $[{\g}_1,{\g}_2]$ be the subspace spanned by $[u,v]$, where $u\in{\g}_1$
and $v\in{\g}_2$. Then
\begin{equation} \label{ga}
[\mathcal{S},\mathcal{S}] \subset \mathcal{A}, \quad [\mathcal{A},\mathcal{A}]\subset\mathcal{A}, \quad [\mathcal{A},\mathcal{S}]=[\mathcal{S},\mathcal{A}]\subset\mathcal{S}.
\end{equation}
According to \eqref{Eq_covder}, we have
\begin{equation} \label{uv}
\nabla_uv = \begin{cases}
\frac{1}{2}[u, v], & \text{if } u, v \in \mathcal S \text{ or } u, v \in \mathcal{A}\\
-\frac{1}{2}[u, v], & \text{if } u \in \mathcal S, v \in \mathcal A\\
\frac{3}{2}[u, v], & \text{if } u \in \mathcal A, v \in \mathcal S.
\end{cases}
\end{equation}

\begin{lemma}\label{lma_w}  Given $u, v, w\in \mathfrak {g}$, we have
\[
\langle [w,u],v\rangle = \langle u,  [w^\top, v]\rangle.
\]
In particular,
\[
\langle [w, u], v \rangle = \begin{cases}
\langle u, [w, v]\rangle, & \text{if } w \in \mathcal{S}\\
- \langle u, [w, v]\rangle, & \text{if } w \in \mathcal{A}.
\end{cases}
\]
\end{lemma}
\begin{proof}
Given $u, v, w\in \mathfrak {g}$, 
\[
\langle [w, u], v \rangle = \mathrm{tr}[(wu-uw) ^\top v] = \mathrm{tr}(u^\top w^\top v) - \mathrm{tr}(w^\top u^\top v) = \mathrm{tr} (u^\top w^\top v) - \mathrm{tr}(u^\top vw^\top) = \langle u,[w^\top,v]\rangle.
\]
\end{proof}


\section{Sectional curvature}\label{sectional}

As stated in the previous section, let $G=\GL(n,\R)$ be equipped with the left invariant Riemannian metric determined by the Frobenius inner product on the Lie algebra ${\g}$ of $G$. The curvature tensor
of the Riemannian connection is given by \cite[p.43]{Helgason}
\begin{equation} \label{curvten}
R(u, v)w = \nabla_u(\nabla_v w) -  \nabla_v(\nabla_u w) - \nabla_{[u, v]} w
\end{equation}
for $u,v,w\in{\g}$. The curvature tensor is to measure the intrinsic bending of $G$ and the bending at each point is measured by the failure of mixed partial derivatives to commute. The sectional curvature of the section spanned by linearly independent $u$ and $v$ in ${\g}$ is
\begin{equation} \label{Suv}
S(u,v) = \frac{\langle R(u,v)v,u\rangle}{\langle u,u\rangle\langle v,v\rangle - \langle u,v\rangle^2}.
\end{equation}
Note that the denominator is always positive. Indeed it is the area $|u\wedge v|$ of the parallelogram determined by $u$ and $v$, so the sign of $S(u,v)$ is the same as that of $\langle R(u,v)v,u\rangle$. Moreover,
when $u$ and $v$ are orthonormal, $S(u,v)=\langle R(u,v)v,u\rangle$.

When $u$ and $v$ are regarded as left invariant vector fields, the sectional curvature $S(u,v)$ is in general a function on $G$, but by the left invariance of the metric, it is a constant on $G$.

We will compute $\langle R(u,v)v,u\rangle$ for any $u, v \in {\g}$, starting with some special cases. Recall that $ \mathcal{S}$ and $ \mathcal{A}$ are the spaces of symmetric
and skew-symmetric matrices in ${\g}$, respectively.

\begin{theorem} \label{Th_uv_sym}
	Let $u, v\in \mathcal{S}$. Then $\langle R(u,v)v,u\rangle = - \frac 74\norm{[u, v]}^2 \leq  0$.
\end{theorem}

\begin{proof}
Since $u,v\in\mathcal{S}$ and $[u,v]\in \mathcal{A}$, we have, by (\ref{uv}) and (\ref{curvten}),
\[\begin{aligned}
R(u, v)v &= \nabla_u(\nabla_v v) - \nabla_v(\nabla_u v) -\nabla_{[u, v]} v\\
&= 0 + \frac{1}{4}[v, [u,v]]- \frac{3}{2}[[u, v], v]\\
& = - \frac{7}{4}[[u, v], v].
\end{aligned}
\]
Then by Lemma~\ref{lma_w},
\[
\langle R(u, v)v, u\rangle = \langle - \frac{7}{4}[[u, v], v], u \rangle = -\frac{7}{4} \langle [u, v], [u,v] \rangle = - \frac{7}{4} \norm{[u, v]}^2 \leq  0.
\]
\end{proof}

\begin{theorem} \label{Th_uv_skew}
	Let $u, v\in   \mathcal{A}$. Then $\langle R(u,v)v,u\rangle = \frac 14 \norm{[u, v]}^2 \geq  0$.
\end{theorem}

\begin{proof}
This is proved in the same way as for Theorem~\ref{Th_uv_sym}. We first obtain $R(u,v)v=-\frac 14[[u,v],v]$ and then we get
 $\langle R(u,v)v,u\rangle=\frac 14\langle [u,v],[u,v]\rangle$.
\end{proof}

\begin{theorem} \label{Th_uv_symskew}
Let $u \in  \mathcal{S}$ and $v \in \mathcal{A}$. Then $\langle R(u,v)v,u\rangle =  \frac 14 \norm{[u, v]}^2 \geq  0$.
\end{theorem}

\begin{proof}
This is proved in the same way as for Theorem~\ref{Th_uv_skew}.
\end{proof}

\begin{claim} \label{Clm_X1X2_zero}
Let $u_1 \in  \mathcal{S}$ and $u_2 \in \mathcal{A}$. Then $\langle R(u_1, v)v, u_2 \rangle = \langle R(u_2, v)v, u_1 \rangle = 0$ for $v\in \mathcal{A}\cup \mathcal{S}$.
\end{claim}

\begin{proof}
Assume $v\in \mathcal{A}$. Because $[u_1, v]\in \mathcal{S}$  and $[u_2, v]\in \mathcal{A}$, $[u_1, v]$ and $[u_2, v]$ are orthogonal, i.e., $\langle [u_1, v], [u_2,v] \rangle = 0$.
Then we have
\[
\langle R(u_1, v)v, u_2 \rangle =  \langle - \frac{1}{4}[[u_1, v], v], u_2 \rangle = \frac{1}{4} \langle [u_1, v], [u_2,v] \rangle = 0.
\]
and
\[
\langle R(u_2, v)v, u_1 \rangle = \langle - \frac{1}{4}[[u_2, v], v], u_1 \rangle = \frac{1}{4} \langle [u_2, v], [u_1,v] \rangle = 0.
\]
The proof for $v\in \mathcal{S}$ is similar.
\end{proof}

For any $u \in {\g}$, we will let $u_1=(u+u^\top )/2$ and $u_2=(u-u^\top )/2$. Then $u=u_1+u_2$ is the decomposition ${\g}=\mathcal{S}\oplus\mathcal{A}$.

\begin{theorem} \label{Th_uv_arbskew}
Let $u \in {\g}$ and $v \in \mathcal{A}$. Then $\langle R(u, v)v, u\rangle = \frac 14\norm{[u,v]}^2\geq 0$.
\end{theorem}

\begin{proof}
For any $u\in{\g}$, we have
\[
\begin{aligned}
\langle R(u, v)v, u\rangle & = \langle R(u_1+u_2, v)v, u_1+u_2 \rangle \\
& = \langle R(u_1, v)v, u_1 \rangle + \langle R(u_2, v)v, u_2 \rangle + \langle R(u_1, v)v, u_2 \rangle + \langle R(u_2, v)v, u_1 \rangle \\
& = \frac{1}{4} \norm{[u_1, v]}^2 +\frac{1}{4} \norm{[u_2, v]}^2 \quad \mbox{(by  Theorems \ref{Th_uv_skew} and \ref{Th_uv_symskew}, and Claim~\ref{Clm_X1X2_zero})} \\
& = \frac{1}{4}\norm{[u,v]}^2  \geq 0.
\end{aligned}
\]
\end{proof}

\begin{theorem} \label{Th_uv_arbsym}
Let $u \in {\g}$ and $v \in \mathcal{S}$. Then
\[
\langle R(u, v)v, u\rangle = - \frac{7}{4} \norm{[u_1, v]}^2 +\frac{1}{4} \norm{[u_2, v]}^2.
\]
\end{theorem}

\begin{proof}
For any $u \in {\g}$, we have
\[
\begin{aligned}
\langle R(u, v)v, u \rangle &= \langle R(u_1+u_2, v)v, u_1+u_2 \rangle \\
&= \langle R(u_1, v)v, u_1 \rangle + \langle R(u_2, v)v, u_2 \rangle + \langle R(u_1, v)v, u_2 \rangle + \langle R(u_2, v)v, u_1 \rangle \\
&= - \frac{7}{4} \norm{[u_1, v]}^2 +\frac{1}{4} \norm{[u_2, v]}^2 \quad \mbox{(by Theorems \ref{Th_uv_sym} and \ref{Th_uv_symskew}, and Claim~\ref{Clm_X1X2_zero}).}
\end{aligned}
\]
Note that we have used $\langle R(u_2,v)v,u_2\rangle=\langle R(v,u_2)u_2,v\rangle$ from the standard curvature identity \cite [p.69]{Helgason}:
\begin{equation} \label{curvprop}
\langle R(X,Y)Z,W\rangle = \langle R(Y,X)W,Z\rangle.
\end{equation}
\end{proof}


To obtain $\langle R(u,v)v,u\rangle$ for any $u,v\in{\g}$, we need to prove the following claim first.

\begin{claim}\label{Clm_XY_equation}
For any $u, v \in {\g}$, we have
\[
\norm{[u, v]}^2 = \norm{[u, v_1]}^2+\norm{[u, v_2]}^2-2\langle [v_1, v_2], [u_1, u_2]\rangle.
\]
\end{claim}

\begin{proof}
\[
	\begin{aligned}
	\norm{[u, v]}^2 
	&=\langle [u, v_1+v_2],[u, v_1+v_2]\rangle\\
	&=\langle [u, v_1], [u, v_1]\rangle+\langle [u, v_2], [u, v_2]\rangle +2\langle [u, v_1], [u, v_2]\rangle\\
	&= \norm{[u, v_1]}^2+\norm{[u, v_2]}^2+2\{\langle [u_1, v_1], [u_2, v_2]\rangle+\langle [u_2, v_1], [u_1, v_2]\rangle\}\\
	& \mbox{(because $[u_1,v_1]$ and $[u_1,v_2]$ are orthgonal, and so are $[u_2,v_1]$ and $[u_2,v_2]$)} \\
	&=\norm{[u, v_1]}^2+\norm{[u, v_2]}^2+2\{\langle-[u_2 ,[u_1, v_1]], v_2\rangle+\langle [u_1, [u_2, v_1]],v_2\rangle\}\\
	&=\norm{[u, v_1]}^2+\norm{[u, v_2]}^2+2\{\langle[u_2 ,[v_1, u_1]], v_2\rangle+\langle [u_1, [u_2, v_1]],v_2\rangle\}\\
	&=\norm{[u, v_1]}^2+\norm{[u, v_2]}^2+2\langle-[v_1,[u_1, u_2]], v_2\rangle \quad \mbox{(by Jacobi identity~\eqref{Eq_Jacobi})}\\
	&=\norm{[u, v_1]}^2+\norm{[u, v_2]}^2-2\langle[u_1, u_2],[v_1, v_2]\rangle.
	\end{aligned}
\]
\end{proof}

\begin{theorem} \label{Th_general}
	Let $u, v \in {\g}$. Then
	\begin{equation} \label{Eq_XY_general}
	\langle R(u, v)v, u \rangle =  - 2 \norm{[u_1, v_1]}^2 +\frac{1}{4} \norm{[u, v]}^2 +2 \langle [u_1, v_1], [u_2, v_2]\rangle.
	\end{equation}
\end{theorem}

\begin{proof}
By (\ref{curvprop}), we have
\begin{equation} \label{Ruvvu}
	\begin{aligned}
		\langle R(u, v)v, u \rangle & = \langle R(u, v_1+v_2)(v_1+v_2), u \rangle \\
		& = \langle R(u, v_1)v_1, u \rangle+\langle R(u, v_2)v_2, u \rangle+2\langle R(u, v_1)v_2, u \rangle.
	\end{aligned}
\end{equation}
The last term above without factor $2$ is
\[
	\begin{aligned}
		\langle R(u, v_1)v_2, u \rangle & = \langle R(u_1+u_2, v_1)v_2, u_1+u_2 \rangle \\
		& = \langle R(u_1, v_1)v_2, u_1 \rangle+\langle R(u_2, v_1)v_2, u_2 \rangle+\langle R(u_1, v_1)v_2, u_2 \rangle +\langle R(u_2, v_1)v_2, u_1 \rangle.
	\end{aligned}
\]
To compute the second and fourth terms, we first find $ R(u_2, v_1)v_2$ using (\ref{curvten}) and (\ref{uv}).
\[
	\begin{aligned}
		R(u_2, v_1)v_2 & =  \nabla_{u_2}(\nabla_{v_1} v_2) -  \nabla_{v_1}(\nabla_{u_2} v_2) - \nabla_{[u_2, v_1]} v_2\\
		& = -\frac{3}{4}[u_2,[v_1, v_2]]-\frac{1}{4}[v_1, [v_2, u_2]]+\frac{1}{2}[v_2,[v_1,u_2]].
	\end{aligned}
\]
Since $[v_1, u_2]$ and $[v_2, u_2]$ are orthogonal, we have
\[
\langle R(u_2, v_1)v_2, u_2 \rangle = \frac{3}{4}\langle [v_1,v_2],[u_2,u_2]\rangle - \frac{1}{4}\langle [v_2,u_2],[v_1,u_2]\rangle - \frac{1}{2} \langle [v_1, u_2], [v_2, u_2]\rangle = 0,
\]
and
\[
	\begin{aligned}
		  & \langle R(u_2, v_1)v_2, u_1 \rangle \\
		  = &\,\,-\frac{3}{4}\langle [u_2,[v_1, v_2]], u_1\rangle - \frac{1}{4}\langle [v_1, [v_2, u_2]], u_1\rangle+\frac{1}{2}\langle[v_2, [v_1, u_2]], u_1\rangle\\
		  = &\,\,-\frac{1}{4}\langle [u_2,[v_1, v_2]], u_1\rangle-(\frac{1}{2}\langle [u_2,[v_1, v_2]], u_1\rangle+\frac{1}{2}\langle[v_2, [u_2, v_1]], u_1\rangle )- \frac{1}{4}\langle [v_1, [v_2, u_2]], u_1\rangle\\
		  = &\,\,-\frac{1}{4}\langle [u_2,[v_1, v_2]], u_1\rangle +\frac{1}{2}\langle v_1, [v_2, u_2]], u_1\rangle- \frac{1}{4}\langle [v_1, [v_2, u_2]], u_1\rangle \quad \mbox{(by Jacobi identity~\eqref{Eq_Jacobi})} \\
		   = &\,\,-\frac{1}{4}\langle [v_1, v_2], [u_1, u_2]\rangle+\frac{1}{4} \langle [v_1, u_1], [v_2, u_2]\rangle.
	\end{aligned}
\]
The computation of the first and third terms in $\langle R(u,v_1)v_2,u\rangle$ is similar. We first find
\[
	\begin{aligned}
		R(u_1, v_1)v_2 &=  \nabla_{u_1}(\nabla_{v_1} v_2) -  \nabla_{v_1}(\nabla_{u_1} v_2) - \nabla_{[u_1, v_1]} v_2\\
		& = -\frac{1}{4}[u_1,[v_1, v_2]]+\frac{1}{4}[v_1, [u_1, v_2]]-\frac{1}{2}[[u_1,v_1],v_2].
	\end{aligned}
\]
Since $[u_1, v_1]$ and $[u_1, v_2]$ are orthogonal, we have
\[
 \langle R(u_1, v_1)v_2, u_1 \rangle = -\frac{1}{4}\langle [v_1,v_2],[u_1,u_1]\rangle + \frac{1}{4}\langle [u_1,v_2],[v_1,u_1]\rangle - \frac{1}{2} \langle [u_1, v_1], [u_1,v_2]\rangle = 0,
\]
and
\[
	\begin{aligned}
		  &\langle R(u_1, v_1)v_2, u_2 \rangle \\
		  =  &\,\, -\frac{1}{4}\langle [u_1,[v_1, v_2]], u_2\rangle + \frac{1}{4}\langle [v_1, [u_1, v_2]], u_2\rangle - \frac{1}{2}\langle [[u_1,v_1],v_2], u_2 \rangle\\
		  = & \,\, \frac{1}{4}\langle [u_1,[v_2, v_1]], u_2\rangle+\frac{1}{4}\langle [v_1, [u_1, v_2]], u_2\rangle - \frac{1}{2}\langle [[u_1,v_1],v_2], u_2 \rangle\\
		  = & \,\, -\frac{1}{4}\langle [v_2,[v_1, u_1]], u_2\rangle -\frac{1}{2}\langle [[u_1,v_1],v_2], u_2 \rangle \quad \mbox{(by Jacobi identity~\eqref{Eq_Jacobi})} \\
		  = & \,\, \frac{1}{4}\langle [v_1, u_1], [v_2, u_2]\rangle +\frac{1}{2}\langle [u_1,v_1], [u_2,v_2]\rangle\\
		 = & \,\, \frac{3}{4}\langle [v_1, u_1], [v_2, u_2]\rangle.
	\end{aligned}
\]
So, we obtain
\[
	\begin{aligned}
		&\langle R(u, v_1)v_2, u \rangle \\
		= & \,\, \frac{3}{4} \langle [v_1, u_1], [v_2, u_2]\rangle - \frac{1}{4} \langle [v_1, v_2], [u_1, u_2]\rangle+\frac{1}{4} \langle [v_1, u_1], [v_2, u_2]\rangle \\
		= & \, \,\langle [v_1, u_1], [v_2, u_2]\rangle  - \frac{1}{4} \langle [v_1, v_2], [u_1, u_2]\rangle.
	\end{aligned}
\]
We now return to (\ref{Ruvvu}). Applying Theorem~\ref{Th_uv_arbsym}  to $\langle R(u,v_1)v_1,u\rangle$ and applying Theorem~\ref{Th_uv_arbskew}  to $\langle R(u,v_2)v_2,u\rangle$,
we have
\[
	\begin{aligned}
		& \langle R(u, v)v, u \rangle \\
		= &\,\, - \frac{7}{4} \norm{[u_1, v_1]}^2 +\frac{1}{4} \norm{[u_2, v_1]}^2 +\frac{1}{4} \norm{[u, v_2]}^2 +2 \langle [v_1, u_1], [v_2, u_2]\rangle-\frac{1}{2} \langle [v_1, v_2], [u_1, u_2]\rangle\\
		= &\,\, - 2 \norm{[u_1, v_1]}^2 +\frac{1}{4} \norm{[u, v_1]}^2 +\frac{1}{4} \norm{[u, v_2]}^2 +2 \langle [v_1, u_1], [v_2, u_2]\rangle-\frac{1}{2} \langle [v_1, v_2], [u_1, u_2]\rangle\\
		& \mbox{(because $\norm{[u_1,v_1]}^2+\norm{[u_2,v_1]}^2=\norm{[u,v_1]}^2$)} \\
		= &\,\, - 2 \norm{[u_1, v_1]}^2+\frac{1}{4}\{\norm{[u, v_1]}^2+\norm{[u, v_2]}^2-2\langle[u_1, u_2],[v_1, v_2]\rangle\}+2\langle [v_1, u_1], [v_2, u_2]\rangle\\
		= &\,\, - 2 \norm{[u_1, v_1]}^2 +\frac{1}{4} \norm{[u, v]}^2 +2 \langle [v_1, u_1], [v_2, u_2]\rangle \quad \mbox{(by Claim~\ref{Clm_XY_equation})}.
	\end{aligned}
\]
\end{proof}

Let $H$ be a closed subgroup of $G=\GL(n,\R)$. Then $H$ is a Lie subgroup of $G$ and its Lie algebra ${\h}$ is a sub-Lie algebra of ${\g}$. The left invariant Riemannian metric on $G$ inducesa left invariant Riemannian metric on $H$
by restricting to the tangent spaces of $H$. Then $H$ becomes a sub-Riemannian manifold of $G$. For any $u,v\in{\h}$,
we may compute the sectional curvature $S(u,v)$ on $G$ as defined by (\ref{Suv}), and we may also compute the sectional curvature $S_H(u,v)$ on $H$. In general, they are different, but if $H$ is a totally geodesic sub-manifold
of $G$, then $S(u,v)=S_H(u,v)$. By definition, $H$ is a totally geodesic sub-manifold of $G$ if all the geodesics in $G$, starting in $H$ and tangent to $H$, are contained in $H$ and so are also
geodesics in $H$. In this case, it is well known that $S(u,v)=S_H(u,v)$; see for example Theorem~12.2 in \cite[Chapter~I]{Helgason}.

By (3.9) in \cite{Martin2016}, the geodesic $\gamma(t)$ in $G=\GL(n,\R)$ with $\gamma(0)=e$ and $\gamma'(0)=u\in{\g}$ is given by
\[\gamma(t) = \exp(tu^\top )\exp(t(u-u^\top )).\]
Assume $H$ is transpose-invariant, that is, for any $h\in H$, $h^\top \in H$. Then ${\h}$ is also transpose-invariant. From the above geodesic expression,
it is clear that any geodesic in $G$ emitting from $e$ and tangent to $H$ is contained in $H$. Because the Riemannian metric is left invariant, this is true for the geodesic emitting from any point in $H$.
It follows that $H$ is a total geodesic sub-manifold of $G$. We have proved the following result.

\begin{theorem} \label{thH}
Let $H$ be a closed and transpose-invariant subgroup of $G=\GL(n,\R)$, and let it be equipped with the left invariant Riemannian metric determined by the Frobenius inner product restricted to its Lie algebra ${\h}$. Then
Theorem~\ref{Th_general} holds on $H$, that is, (\ref{Eq_XY_general}) holds for $u,v\in{\h}$.
\end{theorem}

\section{Zero curvature and commuting matrices in ${\gl}(n,\R)$}
In this section let ${\g} = {\gl}(n,\R)$. 
Let $u$ and $v$ be two commuting matrices in ${\g}$. Let $A$ be the abelian Lie subgroup of $G = \GL (n, \R)$ with Lie algebra span$\{u,v\}$, and let $A$ be equipped
with the induced Riemannian metric from $G$. Because $A$ is abelian, by (\ref{covderiv}), the covariant derivative of invariant vector fields vanishes on $A$,
and  hence it has a zero curvature. It is interesting to know whether the sectional curvature $S(u,v)$ on $G$ of the section spanned by $u$ and $v$ is also zero.
It is also interesting to know, if $S(u,v)=0$, whether $u$ and $v$ commute.
We will see that the answers to both questions are negative in general, but are positive under additional conditions.

First we note from the following theorem that $S(u,v)\leq 0$ if $u$ and $v$ commute.

\begin{theorem}\label{Th_XY_commute}
If $[u,v]=0$, then $\langle R(u,v)v,u\rangle = -4\norm{[u_1,v_1]}^2 \leq 0$.
\end{theorem}

\begin{proof}
We have
\[
	[u, v] = [u_1+u_2, v_1+v_2] = [u_1, v_1]+[u_2, v_2]+[u_1, v_2]+[u_2, v_1].
\]
Since $[u, v] = 0$, its skew-symmetric part $[u_1, v_1]+[u_2, v_2]$ must be zero. We have $[u_1, v_1]= - [u_2, v_2]$. By Theorem~\ref{Th_general},
\[\langle R(u,v)v,u\rangle = -2\norm{[u_1,v_1]}^2 + 0 - 2\langle [u_1,v_1],[u_1,v_1]\rangle = -4\norm{[u_1,v_1]}^2.\]
\end{proof}

The following theorem is a direct consequence of Theorems \ref{Th_uv_sym} and \ref{Th_uv_skew}.

\begin{theorem}
Let $u,v\in{\g}$. Assume either of the following two conditions:
\begin{itemize}
\item[(i)] both $u$ and $v$ are symmetric; or
\item[(ii)]  either $u$ or $v$ is skew-symmetric.
\end{itemize}
Then $[u,v]=0$ if and only if $S(u,v)=0$.
\end{theorem}

In general, there are non-commuting matrices $u$ and $v$ with zero sectional curvature $S(u,v)$, such as
\[
u=
\begin{bmatrix} 
1&\sqrt{7}/2\\
-\sqrt{7}/2 & 2
\end{bmatrix}
 \quad {\rm and} \quad
v=
\begin{bmatrix}
0&1\\
1&0
\end{bmatrix},
\]
and there are commuting matrices $u$ and $v$ with negative sectional curvature $S(u,v)$, such as
\[u = \left[\begin{array}{rrr} 1 & 1 & -1 \\ 1 & 1 & 0 \\ 2 & 0 & 1 \end{array} \right] \quad {\rm and} \quad
v = \left[\begin{array}{rrr} 0 & -1 & 1 \\ -1 & 2 & -1 \\ -2 & 2 & -1 \end{array} \right].\]
We note that it is not possible to have $2\times 2$ commuting matrices $u$ and $v$ with a nonzero sectional curvature. To see this, note that any two $2\times 2$ skew-symmetric matrices commute. Then by the proof
of Theorem~\ref{Th_XY_commute}, $[u_1,v_1]=-[u_2,v_2]=0$, which implies $S(u,v)=0$ by Theorem~\ref{Th_XY_commute}.

\section{Reductive Lie group}

Let us recall the definition of reductive group \cite [Chapter VII] {Knapp}.

\begin{definition} The {\it Harish-Chandra class} $\mathcal H$ consists of $4$-tuples $(G, K, \theta, B)$, where $G$ is a Lie group, $K$ is a compact subgroup of $G$, $\theta$ is a Lie algebra involution of the Lie algebra $\g$ of $G$, and $B$ is a nondegenerate, $\Ad (G)$-invariant, symmetric, bilinear form on $\g$ such that
\begin{enumerate}
\item $\g$ is reductive, i.e., $\g=\g_1+\z$, where $\g_1=[\g,\g]$ and $\z$ is the center of $\g$.
\item $\g=\k+\p$ (called the {\it{Cartan decomposition}}), where $\k$ is the $+1$-eigenspace and $\p$ is the $-1$-eigenspace under  $\theta$.
\item $\k$ and $\p$ are orthogonal with respect to $B$, and $B$ is negative definite on $\k$ and positive definite on $\p$.
\item the map $K\times\exp\p\to G$ given by multiplication is a surjective diffeomorphism.
\item for every $g\in G$, the automorphism $\Ad (g)$ of $\g$, extended to the complexification $\g^{\C}$ of $\g$ is contained in $\Int \g^{\C}$.
\item the analytic subgroup $G_1$ of $G$ with Lie algebra $\g_1=[\g, \g]$ has finite center.
\end{enumerate}
\end{definition}
If $(G, K, \theta, B)\in\mathcal H$, then $G$ is called a {\it reductive Lie group}.

The bilinear form $B(\cdot,
\cdot):\g\times \g\to \R$
  induces an Euclidean
inner product $B_{\theta}(\cdot, \cdot)$ on $\g$ \cite [p.448]{Knapp}:
\[
\langle X,Y\rangle :=B_{\theta}(X, Y) = -B(X, \theta Y).
\]
Note that $
B_{\theta}|_{(\k \times \k)} = -B \quad \text{and} \quad B_{\theta}|_{(\p \times \p)} = B,
$
and that $\k$ and $\p$ are orthogonal under $B$ and thus under $B_{\theta}$ \cite {Knapp}. 

\begin{example} \label{ex} $G=\GL(n,\R)$ ($\GL(n,\C)$) is reductive with $B(X,Y)=\tr (XY)$ ($B(X,Y)=\mathrm{Re}[\tr (XY)]$, where $\mathrm{Re}$ is the real part) and  $\theta X= -X^\top$ ($\theta X = -X^*$). Then $\langle X, Y \rangle =
B_{\theta}(X,Y)=  \mathrm{tr}(X^\top Y)$, where $X, Y\in \gl(n,\R)$ ($\langle X, Y \rangle =\mathrm{Re}[\mathrm{tr}(X^* Y)]$ for $X,Y\in\gl(n,\C)$).
\end{example}

As in Section 3, the covariant derivative $\nabla_uv$ under a left invariant metric is  \cite[(5.3)]{Milnor} 
\begin{equation} \label{covderiv2}
\langle\nabla_uv,w\rangle = \frac{1}{2}(\langle [u,v],w\rangle - \langle [v,w],u\rangle - \langle [u,w],v\rangle).
\end{equation}

In the rest of the paper, we will assume that $G$ is equipped with the left invariant Riemannian metric determined by $\langle\cdot,\cdot\rangle=B_\theta(\cdot, \cdot)$ at ${\g}$.
Let $\norm{u}=\langle u,u\rangle^{1/2}$ be the associated norm.
It is easy to show 
\begin{equation}\label{theta2}
\langle [u,w],v\rangle= - \langle w,[\theta u,v]\rangle.
\end{equation}

\begin{proof}
$
\langle [u, w], v \rangle = -B([u,w],\theta v) 
=- B(w, [\theta v, u])  
= -B(w, \theta [v, \theta u])
= B_\theta(w, [v, \theta u])
= - \langle w,[\theta u,v]\rangle.
$
\end{proof}

By (\ref{covderiv2}) and (\ref{theta2}),
\[\langle\nabla_uv,w\rangle = \frac{1}{2}(\langle [u,v],w\rangle + \langle [\theta v,u],w\rangle + \langle [\theta u,v],w\rangle).\]
It follows that
\begin{equation}\label{Eq_covder2}
	\nabla_uv = \frac{1}{2}([u,v] -[u,\theta v] -[v,\theta u]).
\end{equation}

It is easy to show that
\begin{equation} \label{ga}
[\k ,\k ]\subset\k  , \quad  [\p ,\p ] \subset \k , \quad [\k ,\p ]=[\p ,\k ]\subset\p .
\end{equation}
According to \eqref{Eq_covder2}, we have
\begin{equation} \label{uv2}
\nabla_uv = \begin{cases}
\frac{1}{2}[u, v], & \text{if } u, v \in \p  \text{ or } u, v \in \k \\
-\frac{1}{2}[u, v], & \text{if } u \in \p , v \in \k\\
\frac{3}{2}[u, v], & \text{if } u \in \k, v \in \p .
\end{cases}
\end{equation}

By \eqref{theta2}, $\langle [w,u],v\rangle = - \langle u,  [\theta w, v]\rangle$ and we have the following lemma.
\begin{lemma}\label{lma_w}  Given $u, v, w\in \mathfrak {g}$, we have
\[
\langle [w, u], v \rangle = \begin{cases}
 \langle u, [w, v]\rangle, & \text{if } w \in \p \\
-\langle u, [w, v]\rangle, & \text{if } w \in \k .
\end{cases}
\]
\end{lemma}

The curvature tensor $R$
and sectional curvature  $S$ are defined in the same ways as
\eqref{curvten} and \eqref{Suv}, respectively.
We have the following result and we skip the proofs which are similar to those in Section \ref{sectional}.

\begin{theorem} \label{Th_general_reductive} Let $G$ be a reductive Lie group.
Let $u, v \in {\g}$. Then
	\begin{equation} \label{Eq_XY_general}
	\langle R(u, v)v, u \rangle =  - 2 \norm{[u_1, v_1]}^2 +\frac{1}{4} \norm{[u, v]}^2 +2 \langle [u_1, v_1], [u_2, v_2]\rangle.
	\end{equation}

So 
\begin{enumerate}
\item Let $u, v\in \p $. Then $\langle R(u,v)v,u\rangle = - \frac 74\norm{[u, v]}^2 \leq  0$.
\item Let $u, v\in   \k $. Then $\langle R(u,v)v,u\rangle = \frac 14 \norm{[u, v]}^2 \geq  0$.
\item Let $u \in  \p $ and $v \in \k $. Then $\langle R(u,v)v,u\rangle =  \frac 14 \norm{[u, v]}^2 \geq  0$.
\item Let $u \in {\g}$ and $v \in \k $. Then $\langle R(u, v)v, u\rangle = \frac 14\norm{[u,v]}^2\geq 0$.
\item Let $u \in {\g}$ and $v \in \p $. Then
\[
\langle R(u, v)v, u\rangle = - \frac{7}{4} \norm{[u_1, v]}^2 +\frac{1}{4} \norm{[u_2, v]}^2.
\]
\end{enumerate}
\end{theorem}

Let $H$ be a closed subgroup of  $G=\GL(n,\R)$ (or $G=\GL(n,\C)$) that is invariant under (conjugate) transposition. It is known that \cite [p.447]{Knapp}  $H$ is a reductive Lie group. By Theorem~\ref{Th_general_reductive},
we obtain an alternative proof of Theorem~\ref{thH} and extend it to include complex matrix groups.

\begin{cor} \label{thH_reductive}
Let $H$ be a closed  subgroup of $\GL(n,\R)$ or ($\GL(n,\C)$) that is invariant under (conjugate) transposition, and let it be equipped with the left invariant Riemannian metric determined by the inner product
in Example~\ref{ex} restricted to its Lie algebra. Then Theorem~\ref{Th_general} holds for $H$.
\end{cor}

\begin{remark}\rm Let $G$ be a reductive Lie group with Lie algebra $\g$. Let $\g=\k+\p$ be a given Cartan decomposition corresponding to the Cartan involution $\theta$. Let $K$ be the analytic subgroup of $G$ with Lie algebra $\k$. Let $P=\exp \p$. Note that $G=K\exp \p$ and $P$ is not a group, so Theorem~\ref{Th_general_reductive} does not apply. 
When $P$ is equipped with the symmetric space metric, it is a Riemannian manifold and the geodesic starting from $p\in P$ takes the form $p^{1/2} \exp(tu) p^{1/2}$, $u\in \p$ \cite {Ming}. It is related to the geometric means in the context of symmetric space of noncompact type and \cite {Ming} evolves from the study of the matrix geometric means of two $n\times n$ positive definite matrices \cite {Bhatia}. See \cite {Lawson2013, Lawson2014} for some recent interesting results and generalizations of matrix geometric means.
\end{remark}

\bigskip
\bibliographystyle{abbrv}
\bibliography{refs}

\end{document}